\documentclass{amsart}
\usepackage{amssymb}

\usepackage{latexsym}
\usepackage{amsmath}
\usepackage{amssymb}
\usepackage{latexsym}
\usepackage{amsmath}
\usepackage{amscd}
\usepackage{euscript}

\usepackage{graphicx,color}

\usepackage{amsfonts}
\usepackage{times}
\usepackage{bm}

\newtheorem{theorem}{Theorem}[section]

\newtheorem{proposition}[theorem]{Proposition}

\theoremstyle{definition}
\newtheorem{definition}[theorem]{Definition}

\newtheorem{remark}[theorem]{Remark}

\numberwithin{equation}{section}

\begin{document}

\author{C.~Fulton}
\address{Department of Mathematical Sciences \\
Florida Institute of Technology \\
150 W. University Blvd. \\
Melbourne, FL. 32904, USA}
\email{\tt cfulton@fit.edu}
\author{H.~Langer}
\address{Institute for Analysis and Scientific Computing\\
Vienna University of Technology\\
 Wiedner Hauptstrasse 8-10 \\
1040 Vienna, Austria} \email{\tt hlanger@mail.zserv.tuwien.ac.at}
\author{A.~Luger}
\address{Department of Mathematics\\ Stockholm University\\ 106 91 Stockholm, Sweden}

\email{luger@math.su.se}
\thanks{}
\begin{abstract}
\end{abstract}



\title[Directing Functionals]{Mark Krein's Method of Directing Functionals and Singular Potentials
}

\begin{abstract}
It is shown that M.~Krein's method of directing functionals can be used to prove the existence of a scalar spectral measure for certain Sturm--Liouville equations with two singular endpoints. The essential assumption is the existence of a solution of the equation that is square integrable at one singular endpoint and depends analytically on the eigenvalue parameter.
\end{abstract}

\maketitle

\dedicatory{\it Dedicated to Prof.~Eduard Tsekanovskii on the occasion of his 75$^{th}$ birthday}

\section{Introduction}
F.\,Gesztesy and M.\,Zinchenko in \cite{gz}, and, more recently, A.\,Kostenko, A.\,Sakhnovich and G.\,Teschl in  \cite{kst} considered Sturm-Liouville operators $L$ with two singular endpoints which have a {\it scalar} spectral measure. Recall that the classical Weyl theory in the case of two singular endpoints leads to  a $2\times 2$ matrix valued spectral measure and a corresponding $2$-dimensional Fourier transformation. The essential assumption in \cite{gz} and \cite{kst} for the existence of a scalar spectral measure is the existence of a solution of the equation $(L-\lambda)\varphi=0$ which is an  $L^2$--function  for one of the endpoints and depends analytically on the eigenvalue parameter $\lambda$. Then the Fourier transformation to which the scalar spectral measure belongs is defined by means of this solution $\varphi$. The proofs in these papers extend the main ideas of the Titchmarsh-Weyl theory to the present situation; the main difficulty here comes from the fac
 t that
  the Titchmarsh-Weyl function is no longer a Nevanlinna function as in the classical Weyl theory.

It is the aim of the present note to show that, under the assumption of the existence of a holomorphic $L^2$-solution for one endpoint, the existence of a scalar spectral function follows easily from M.G.\,Krein's  method of directing functionals, established more than 60 years ago, see \cite{k}, \cite{k1}. This method also gives a complete answer to the question of the number of spectral functions; combined with M.G.\,Krein's resolvent formula, in the case of nonuniqueness also a description of all these spectral functions can be obtained.

We mention that  I.S.\,Kac in 1956, see \cite{kac1}, was the first to observe that for a Bessel type operator in the case of two singular endpoints  the existence of an analytic $L^2$--solution at one endpoint implies the existence of a scalar spectral function; see also Remark \ref{kacc} at the end of the paper. The authors thank Prof. Kac for calling their attention to his early work in \cite{kac1}, \cite{kac2}, \cite{k3}.

\section{The Method of Directing Functionals}

Let $\mathcal L$ be a complex linear space with a positive semi-definite inner product $(\cdot,\cdot)$. The set $\mathcal D$ is called {\it quasidense} in $\mathcal L$ if for each $f\in\mathcal L$ there exists a sequence $(f_n)$ in $\mathcal D$ such that 
$$
(f-f_n,f-f_n)\longrightarrow 0,\quad n\to\infty.
$$
In the following we consider a linear operator $A_0$ in $\mathcal L$ with quasidense domain ${\rm dom\, A_0}$ that is {\it symmetric}, that is
$$
(A_0f,g)=(f,A_0g),\quad f,g\in{\rm dom}\,A_0.
$$ 

\begin{definition}\label{def1}
A mapping $\Phi$ from $\mathcal L\times \mathbb R$ into the set $\mathbb C$ of complex numbers is called a {\it directing functional} for $A_0$ if it has the following properties:
\begin{itemize}
\item[(1)] 
For each $\lambda\in\mathbb R$, $\Phi(\cdot;\lambda)$ is a linear functional on $\mathcal L$.
\item[(2)]
For each $f\in\mathcal L$, 
$\Phi(f;\cdot)$ is real analytic on $\mathbb R$, that is for each $\lambda_0\in\mathbb R$ it admits a power series expansion in $\lambda-\lambda_0$ with positive radius of convergence.
\item[(3)]
For given $f\in\mathcal L$ and $\lambda\in\mathbb R$ the equation
\begin{equation}\label{1}
A_0g-\lambda g=f
\end{equation}
has a solution $g\in{\rm dom}\,A_0$ if and only if $
\Phi(f;\lambda)=0.$
\end{itemize}
\end{definition}

Clearly, the decisive assumption, that relates the functions $\Phi(f;\cdot)$ with the operator $A_0$, is (3). 
From (1) and (3) it follows easily that the functional $\Phi$ has the property
\begin{equation}\label{2}
\Phi(A_0g;\lambda)=\lambda\,\Phi(g;\lambda),\quad g\in{\rm dom}\,A_0,\ \lambda\in\mathbb R.
\end{equation}

With $\mathcal L_0:=\{f:(f,f)=0\}$, we consider the Hilbert space completion $\mathcal H$ of the factor space $\mathcal L/\mathcal L_0$. The operator $A_0$  induces a densely defined closed symmetric operator $A$ in $\mathcal H$.

 The basic result of M.~Krein can be formulated as follows.

\begin{theorem}\label{t1}
Let $\mathcal L$ be a complex linear space with a positive semi-definite inner product $(\cdot,\cdot)$, and let $A_0$ be  a symmetric linear operator  in $\mathcal L$ with quasidense domain ${\rm dom\, A_0}$ which has  a directing functional $\Phi$. Then:\\
 $(1)$ There exists a measure $\sigma$ on $\mathbb R$ such that the Parseval relation
\begin{equation}\label{5}
(f,f)=\int_\mathbb R\,|\Phi(f;\lambda)|^2\,d\sigma(\lambda),\quad f\in\mathcal L,
\end{equation}
holds.\\
 $(2)$ The measure  $\sigma$ in \eqref{5} is unique if and only if the operator $A$ in $\mathcal H$ is maximal symmetric; otherwise there are infinitely many such measures $\sigma$, which are in a bijective correspondence with all minimal self-adjoint extensions $\widetilde A$ of $A$ in $\mathcal H$ or in a larger subspace $\widetilde {\mathcal H}\supsetneq \mathcal H$.\\
 $(3)$ The mapping $f\mapsto \Phi(f;\cdot)$ extends by continuity to an isometry  from $\mathcal H$ onto a subspace of $L^2_\sigma$; this subspace is the whole of $L^2_\sigma$ if and only if $\widetilde A$ is a self-adjoint extension of $A$  in $\mathcal H$.
\end{theorem}

A self-adjoint extension $\widetilde A$  of a symmetric operator $A$ in $\mathcal H$, which acts in a larger space ${\widetilde{\mathcal H}}\supsetneq \mathcal H$, is called an {\it extension with exit} (from $\mathcal H$).

\begin{remark}\label{uuu} (1) There is no growth restriction of the function $\sigma$ at $\infty$.\\
(2)
The relation \eqref{2} implies that under the mapping $f\mapsto \Phi(f;\cdot)$ from $\mathcal L$ into $L^2_\sigma$ the operator $A$ becomes multiplication by the independent variable.\\
(3) The self-adjoint extension $\widetilde A$ corresponding to the measure $\sigma$ is unitarily equivalent to the operator of multiplication by the independent variable in $L^2_\sigma$. Therefore $\widetilde A$  has simple spectrum.\\
(4) In the original paper \cite{k} M.G.\,Krein considers the more general case of a finite number $n$ of directing functionals (or a mapping $\Phi(\cdot, \lambda)$ from $\mathcal L$ into $\mathbb C^n$). A generalization for a  Hilbert space valued directing functional was given in \cite{l}. In \cite{tex} a situation where the inner product on $\mathcal L$ is indefinite was considered.
\end{remark}

\medskip

For the  proof of Theorem \ref{t1} we refer to \cite{k}, \cite{ag}, \cite{l}.
In applications, $\Phi(f,\lambda)$ is a {\it Fourier transformation for the operator} $A$. The measure  $\sigma$  is called {\it a spectral measure of the operator $A$, corresponding to this Fourier transformation $\Phi$}. 

\section{The Sturm-Liouville Operator}

Consider on the interval $(a,b),\ 0\le a <b\le\infty$, the differential expression $\tau$:
\begin{equation*}
\big(\tau u)(x):=\dfrac{-\left(p(x)u^\prime(x)\right)^\prime+q(x)u(x)}{r(x)},\quad x\in(a,b),
\end{equation*}
with real $q\in L^1_{loc}(a,b)$, $p$ and $r$ continuous and positive in $(a,b)$.
As in \cite{gz} and \cite{kst}, we always assume $\tau$ to be in limit point case at the endpoint $a$, hence there exists an (up to scalar multiples) unique solution of the equation
$$
\tau u=\lambda u \text{ \ on \ } (a,b),\, \Im \lambda\ne 0,
$$
which is square integrable at $a$ (with respect to the weight function $r$). We suppose additionally that this solution depends holomorphically on $\lambda$. More exactly, we make the  following assumption:
 
{\bf(A)}
\textit{There exists a non-trivial solution $u(x,\lambda)=\varphi(x,\lambda)$ of the equation 
$$
\big(\tau u)(x,\lambda)=\lambda\,u(x,\lambda), \qquad x\in(a,b),\,\lambda\in\mathbb C,
$$
with the following properties:\vspace{1mm}
\begin{itemize}
\item[(i)] for each $x\in(a,b)$ the function $\varphi(x,\cdot)$ is entire. \vspace*{1mm}
\item[(ii)] $\varphi(\cdot,\lambda)$ satisfies an $L^2_r$-condition at the endpoint $a$:
$$
\int_a^{b'}|\varphi(x,\lambda)|^2\,r(x)\,dx<\infty
$$
for some $($and hence for all$)$ $a<b'<b$ and all complex $\lambda $.\vspace*{1mm}
\item[(iii)] $\varphi(x,\lambda)$ is real if $x\in(a,b)$ and $\lambda\in\mathbb R$. 
\end{itemize}
}

\begin{remark}
(a) The assumption (iii) is no restriction.\\
(b) In \cite[Hypothesis 3.1]{gz} the $L^2$-property of the solution $\varphi(\cdot,\lambda)$ is assumed to hold only for nonreal $\lambda$. However, in \cite[Examples 3.10 and 3.13]{gz} it holds even for all $\lambda\in\mathbb C$. Moreover, according to \cite[Lemma 2.2]{kst} the apparently weaker assumption in \cite[Hypothesis 3.1]{gz} implies (i) and (ii) above. 

\end{remark}

\medskip

Let $A$ be the minimal closed symmetric operator in $L^2_r(a,b)$ which is associated with $\tau$ in the usual way. Since $\tau$ is limit point at $a$, the operator $A$  has defect $(0,0)$ or $(1,1)$, depending on  $b$ being limit point or limit circle for $\tau$.
Further, $\mathcal L$ denotes the linear manifold of elements of  $L^2_r(a,b)$  for which there is a representative with compact support in $[a,b)$, equipped with the $L^2_r(a,b)$-inner product, and $A_0$ denotes the symmetric operator in $\mathcal L$ given by
\begin{equation}\label{cam}
A_0:=A\Big|_{\mathcal L\cap\,{\rm dom}A}.
\end{equation}

Clearly,  ${\rm dom}\,A_0=\mathcal L\cap\,{\rm dom}\,A$ is dense (and hence quasi-dense) in $\mathcal L$.  

\begin{proposition}\label{main}
If $a$ is limit point for $\tau$ and the assumption $(${\bf A}$)$ is satisfied, then the mapping 
\begin{equation}\label{ddf}
\Phi(f;\lambda):=\int_a^b\varphi(x,\lambda)f(x)\,r(x)\,dx,\quad f\in\mathcal L,\,\lambda\in\mathbb C,
\end{equation}
is a directing functional for $A$.
\end{proposition}

For the convenience of the reader we  note some simple and well-known facts which will be used in the proof of Proposition \ref{main}. Fix any  $c\in (a,b)$ and denote by $\psi_{c,\theta}(\cdot,\lambda)$ the solution of the equation $\tau(u)=\lambda u$ on $(a,c],\ \lambda\in\mathbb C$, which satisfies for $\theta\in[0,\pi)$ the  condition 
$$
\psi_{c,\theta}(c,\lambda)= \sin \,\theta,\quad \psi_{c,\theta}'(c,\lambda) =-\cos\, \theta.
$$


\begin{itemize}
\item[(a)] If $x\in(a,c]$, $\psi_{c,\theta}(x,\cdot)$  is an entire function, satisfying the following condition at $x=c$: 
\begin{equation}\label{bbc}
\cos \theta\, u(c) + \sin\theta\, u'(c)=0,\quad\theta\in[0, \pi ).
\end{equation}

\item[(b)] Under Assumption (A), for $x\in(a,b)$, with $\varphi(x,\cdot)$ also the function $\varphi'(x,\cdot)$ is entire (see \cite{gz}). Indeed, choose any $c\in (a,b)$ and consider the relation
\begin{equation}\label{edi}
\varphi(x,\lambda)=-\varphi'(c,\lambda) \psi_{c,0}(x,\lambda)+ \varphi(c,\lambda)\psi_{c,\pi/2}(x,\lambda),\quad x\in (a,c],\ \lambda\in\mathbb C.
\end{equation}
It implies that $\varphi'(c,\cdot)$ is entire and, differentiating \eqref{edi}, the claim follows.
\item[(c)] Since $\tau$ is limit point at $a$    
the operator $H_{(a,c)}^\theta$, defined in $L^2_r(a,c)$ by $\tau$
and the boundary condition \eqref{bbc}, is self-adjoint.
\item[(d)] If  Assumption ({\bf A}) is satisfied, the resolvent of $H_{(a,c)}^\theta$ is given  by
\begin{equation}\label{res}
\begin{array}{rl}
&\left(\big(H_{(a,c)}^\theta-\lambda\big)^{-1}f\right)(x)\\[2mm]
&\hspace{1.5cm}=\!\dfrac1{w(\lambda)}\!\left[\psi_{c,\theta}(x,\!\lambda)\!\!\displaystyle\int_a^x\!\!\!\!\varphi(\xi,\!\lambda)f(\xi)r(\xi)d\xi\!+\!\varphi(x,\!\lambda)\!\!\displaystyle\int_x^c\!\!\!\!\psi_{c,\theta}(\xi,\!\lambda)f(\xi)r(\xi)d\xi\right]\!\!,
\end{array}
\end{equation}
where $w(\lambda):=\left|\begin{array}{cc}\varphi(x,\lambda)&p(x)\varphi'(x,\lambda)\\\psi(x,\lambda)&p(x)\psi'(x,\lambda)\end{array}\right|$. If we apply this formula to functions $f$ which vanish near the endpoint $a$ the singularites of $\left(\big(H_{(a,c)}^\theta-\lambda\big)^{-1}f\right)$ are the zeros of $w$. Since this set of functions is dense in $L^2_r(a,c)$, the spectrum of $H_{(a,c)}^\theta$ consists of discrete and simple eigenvalues (see also  \cite{kst}), which are strictly decreasing functions of $\theta$. Hence for each $\lambda_0\in\mathbb R$, a $\theta\in[0,\pi)$ can be chosen such that $\lambda_0\in\rho\big(H_{(a,c)}^\theta\big)$.
\end{itemize} 

{\it Proof of Proposition} \ref{main}.
In \eqref{ddf},  $\Phi(f;\lambda)$  is well defined since $f$ vanishes identically near $b$ and, for arbitrary $c\in(a,b)$, $f\big|_{(a,c)},\varphi( \cdot,\lambda)\big|_{(a,c)}\in L^2_r(a,c)$. Clearly, $\Phi(f;\lambda)$ is linear in $f$ and hence (1) in the definition of the directing functional is satisfied. To show (2), we fix $f\in\mathcal L$ and $\lambda_0\in\mathbb R$. 
For $f\in\mathcal L$ we choose $c\in (a,b)$ such that ${\rm supp} f\subseteq[a,c)$, and $\theta\in [0,\pi)$ such that $\lambda_0\in\rho(H_{(a,c)}^\theta)$, see (d). Then the  resolvent of $H_{(a,c)}^\theta$ is holomorphic in $\lambda_0$, and for $x$ with $\max{\rm supp}\, f<x<c$ we obtain 
\begin{equation*}
\left(\big(H_{(a,c)}^\theta-\lambda\big)^{-1}f\right)(x)=\frac{\psi_{c,\theta}(x,\lambda)}{w(\lambda)}\int_a^c\varphi(\xi,\lambda)f(\xi)\,r(\xi)\,d\xi=\frac{\psi_{c,\theta}(x,\lambda)}{w(\lambda)}\Phi(f;\lambda).
\end{equation*}
Since for some $x$ we have $\psi_{c,\theta}(x,\lambda_0)\neq0$  it follows that $\Phi(f;\lambda)$ is analytic at $\lambda_0$.

It remains to show (3). Fix again $\lambda_0\in\mathbb R$ and $f\in\mathcal L$. We have to show that 
$$
\tau g-\lambda_0 g=f\text{ has a solution } g\in{\rm dom}\,A \quad\Longleftrightarrow\quad \Phi(f;\lambda_0)=0.
$$
Assume first that there exists $g\in{\rm dom}\,A$ such that 
\begin{equation}\label{ea}
\tau g-\lambda_0 g=f.
\end{equation}
Since $f,g\in\mathcal L$ we find $c\in (a,b)$ such that ${\rm supp}\,f\subset [a,c)$ and  ${\rm supp}\, g\subset [a,c)$. Choose $\theta\in [0,\pi)$ such that $\lambda_0\in\rho(H_{(a,c)}^\theta)$.  Then \eqref{ea} becomes
 $(H_{(a,c)}^\theta-\lambda_0)g=f$ or $g(x)=\Big(\big(H_{(a,c)}^\theta-\lambda_0\big)^{-1}f\Big)(x)$ on $(a,c]$. For $x$ with $\max{\rm supp}\, g<x<c$ this yields
$$
0=g(x)=\frac{\psi_{c,\theta}(x,\lambda_0)}{w(\lambda_0)}\!\int_a^c\!\!\varphi(\xi,\lambda_0)f(\xi)\,r(\xi)\,d\xi=\frac{\psi_{c,\theta}(x,\lambda_0)}{w(\lambda_0)}\!\int_a^b\!\!\varphi(\xi,\lambda_0)f(\xi)\,r(\xi)\,d\xi,
$$
and hence $\Phi(f;\lambda_0)=0$. 

Conversely, assume  that $\Phi(f;\lambda_0)=0$. We choose again $c$ with ${\rm supp}\, f\subset [a,c)$, and $\theta\in [0,\pi)$ such that $\lambda_0\in\rho(H_{(a,c)}^\theta)$. Then 
$$
g(x):=\left\{\begin{array}{ccc} 
                       \Big(\big(H_{(a,c)}^\theta-\lambda_0\big)^{-1}f\Big)(x) && x<c\\
                       0 && x\geq c
                       \end{array}\right.
$$
belongs to ${\rm dom}\,A$ and is a solution of the equation $\tau g-\lambda_0 g=f$. Indeed, $g\in\mathcal L$, and for $\max\,{\rm supp}\,f<x\le c$ it holds
\begin{align*}
g(x)&= \frac1{w(\lambda_0)}\psi_{c,\theta}(x,\lambda_0)\!\!\int_a^c\!\!\!\varphi(\xi,\lambda_0)f(\xi)\,r(\xi)\,d\xi\\
&=\frac1{w(\lambda_0)}\psi_{c,\theta}(x,\lambda_0)\int_a^b\varphi(\xi,\lambda_0)f(\xi)\,r(\xi)\,d\xi\\
&=\frac1{w(\lambda_0)}\psi_{c,\theta}(x,\lambda_0)\Phi(f;\lambda_0)\\
&=0.
\end{align*}
 Hence  $g$ vanishes identically near $x=c$, therefore $g\in{\rm dom}\,A$ and g satisfies the differential equation.
\qed

An immediate consequence of Proposition \ref{main} and Theorem \ref{t1} is the following theorem.

\begin{theorem}\label{final}
Let $\tau$ be the Sturm-Liouville differential operator
\begin{equation*}
\big(\tau u)(x):=\dfrac{-\left(p(x)u^\prime(x)\right)^\prime+q(x)u(x)}{r(x)},\quad x\in(a,b),
\end{equation*}
with real $q\in L^1_{loc}(a,b)$, $p$ and $r$ continuous and positive in $(a,b)$.
 Suppose that $\tau$ is limt point  at the endpoint $a$, and that the Assumption {\bf(A)} is satisfied. With $\tau$ we associate the symmetric operator $A$ in  $L^2_r(a,b)$ as above and the   Fourier transformation $\Phi(f,\lambda)$ from \eqref{ddf}.  Then for $A$ there exists at least one spectral measure $\sigma$ corresponding to $\Phi(f,\lambda)$. The spectral measure $\sigma$ is unique if and only if $\tau$ is limit point at $b$. If $\tau$ is regular or limit circle at $b$ there are infinitely many spectral measures $\sigma$ for $A$ corresponding to $\Phi(f,\lambda)$; they are in a bijective correspondence with all minimal selfadjoint extensions of the symmetric operator  $A$,  without or with exit from $L^2_r(a,b)$.  
\end{theorem}
 
\begin{remark} (1) To make the last claim more precise, if $A$ is self-adjoint there is exactly one spectral measure $\sigma$. Now suppose that  $A$ has defect $(1,1)$. Then there are infinitely many spectral measures $\sigma$, such that the Fourier transformation maps $L^2_r(a,b)$ onto $L^2(\sigma)$; they can be parametrized through the classical, i.e. $\lambda$--independent self-adjoint boundary conditions at the endpoint $b$. There are also infinitely many spectral measures $\sigma$, such that the Fourier transformation maps $L^2_r(a,b)$ onto a proper subspace of $L^2(\sigma)$;  they can be parametrized through  $\lambda$--dependent self-adjoint boundary conditions at the endpoint $b$, see e.g. \cite{dls}.
\end{remark}

\medskip

As was mentioned already, there is in general no growth restriction for the spectral measure $\sigma$. This is in contrast to the classical situation where a Titchmarsh--Weyl function of Nevanlinna (or Herglotz) class is used which leads to a spectral measure $\sigma$ such that $\displaystyle\int_\mathbb R\,\dfrac{d\,\sigma(\lambda)}{1+\lambda^2}<\infty$. Under more special assumptions information about the growth of $\sigma$ was obtained e.g. in \cite{kac1}, \cite{kac2}, \cite{gz}. In some cases Titchmarsh--Weyl functions were introduced which turned out to be generalized Nevanlinna functions. Then the spectral measure grows polynomially at $\infty$, see e.g. \cite{fl},\cite{kst}.

\medskip

\section{Examples}
As a first example we mention the hydrogen atom  operator
$$
\tau=-\dfrac{d^2}{dx^2}+\left[-\dfrac{a}x+\dfrac{\nu^2-\frac 14}{x^2}\right]\quad 0<x<b,\quad 0<b\le\infty,
$$
in the space $L^2(0,b)$, comp.  \cite{gz}, \cite{fl}, \cite{kl}, \cite{kst},  and \cite{kf}. Assume that $\nu\ge 1$ and $a\in\mathbb R$. Then $\tau$ is limit point at $x=0$, and regular (limit point) at $x=b$ if $0<b<\infty$ ($b=\infty$). A function $\varphi$ satisfying the assumption (A) at $x=0$ is
\begin{equation}\label{erl}
\varphi(x,\lambda)=\dfrac 1{\left(-2{\rm i}\sqrt{\lambda}\right)^{\nu+\frac 12}}\mathcal M_{\beta,\nu}\left(-2{\rm i}x\sqrt{\lambda}\right)=x^{\nu+\frac 12}\left[1+\sum_{n=1}^\infty c_n(\lambda)x^n\right]
\end{equation}
with $\beta:=\dfrac{{\rm i }~a}{2\sqrt \lambda}$, where $c_n(\lambda)$ is the polynomial of degree $\left[\frac n2\right]$ which is defined as the solution of the Frobenius recurrence relation (see \cite[p. 1452]{fu})
$$
c_n(\lambda)=-\dfrac{a}{n(n+2\nu)}c_{n-1}(\lambda)-\dfrac\lambda{n(n+2\nu)}c_{n-2}(\lambda).
$$
 Here
$\mathcal M_{\beta,\nu}$ is the Whittaker function of first kind \cite{sl}.
The directing functional or the Fourier transformation is given by \eqref{ddf} with $\varphi$ from \eqref{erl}.
If $b=\infty$ there is exactly one spectral measure $\sigma$ (it was given explicitly e.g. in \cite{fl}), if $0<b<\infty$ there are infinitely many spectral measures.

\medskip
As a second example we consider the Associated Legendre operator
$$
(\tau y)(x)=-\left((1-x^2)y'(x)\right)' + \left(\dfrac 14+\dfrac{m^2}{1-x^2}\right)y(x),\quad -1<x<1.
$$
If $m\ge 1$ both endpoints $\pm 1$ are limit point.  
A solution $\varphi$ of $\tau y-\lambda y=0$ satisfying the Assumption (A) at $x=-1$ is 
\begin{equation}\label{arth}
\begin{array}{ll}
\varphi(x,\lambda)&\!\!\!=(1-x^2)^{\frac m2}F\left(m+\frac 12-\sqrt{\lambda},m+\frac 12+\sqrt{\lambda},m+1;\dfrac{1+x}2\right)\\
&\!\!\!=(1-x^2)^{\frac m2}\left[1+\sum_{j=1}^\infty\dfrac{\prod_{i=1}^j\left([m+i-\frac 12]^2-\lambda\right)}{j!(m+1)_j}\left(\dfrac{1+x}2\right)^j\right],
\end{array}
\end{equation}
where $(m+1)_j=(m+1)\,m\cdots(m+2-j)$. Here $F$ is the hypergeometric function defined by
$$
F(a,b,c;t)=\sum_{n=0}^\infty\,\dfrac{(a)_n(b)_n}{(c)_n\,n!}\,t^n.
$$
The directing functional is again given by \eqref{ddf} with $\varphi$ from \eqref{arth}, and there is exactly one spectral measure for the problem on $(-1,1)$. 

\medskip

A third example is the Laguerre operator
$$
\tau=x\dfrac{d^2}{dx^2}+(1+\alpha-x)\dfrac d{dx}\qquad\text{on}\ \ (0,b),\quad 0<b\le\infty.
$$
It is symmetric in the space $L^2_r(0,b)$ with weight function $r(x)=x^\alpha e^{-x},\ x\in (0,b)$, and regular at $b$ if $b<\infty$, limit point at $b$ if $b=\infty$. Assume that $\alpha\in\mathbb R,\,|\alpha|> 1$, and $\alpha\ne \pm 2, \pm 3,\dots$. Then $\tau$ is  limit point at zero. The  equation
$$
x\dfrac{d^2y(x)}{dx^2}+(1+\alpha-x)\dfrac {dy(x)}{dx}-\lambda y(x)=0,\quad 0<x<b,
$$
 can also be written as
$$
\left(x^{1+\alpha}e^{-x}y'(x)\right)'-\lambda x^\alpha e^{-x}y(x)=0,\quad x\in (0,b).
$$
A solution $\varphi(x,\lambda)$ satisfying the Assumption (A) at $x=0$ is, for $\alpha\ne -2,-3,\dots$,
$$
\varphi(x,\lambda)=M(\lambda,1+\alpha,x)=1+\sum_{n=1}^\infty\,\dfrac{(\lambda)_n}{(1+\alpha)_nn!}x^n;
$$
here $M$ is the confluent hypergeometric function of first kind. 
 The  directing functional \eqref{ddf} becomes
$$
\Phi(f,\lambda)=\int_0^b\, f(x)\varphi(x,\lambda)\,x^\alpha e^{-x}\,dx.
$$

\medskip

 In all these examples we get spectral measures of polynomial growth at $\infty$, depending on the parameters $\nu,\,m$ or $\alpha$. For the second example and for the third example with $b=\infty$ the corresponding Titchmarsh--Weyl coefficients and spectral functions will be given explicitly in a forthcoming paper. 

\begin{remark}\label{kacc}
If a symmetric operator $A$ has a ($\mathbb C$-valued) directing functional as in Definition \ref{def1}, then the spectrum of each minimal self-adjoint extension of 
the symmetric operator $\widetilde A$ in $\mathcal H$ is simple, also in the case of two singular endpoints, see Remark \ref{uuu},\,(3). The first to observe this for a general class of operators with a Bessel type singularity was I.S.Kac in 1956 \cite{kac1}, see also \cite{kac2}, \cite{k3}. We note also that for the Bessel equation of order $\nu\in(0,\infty)$, Kac's normalization condition in \cite[Theorem 2]{kac1} and \cite[Theorem 13]{kac2},
$$
\lim_{x\to 0},x^{-\nu-\frac 12}\,y(x,\lambda)=1,
$$
selects the same solution $y(x,\lambda)$, entire in $\lambda$, that was recently used in \cite{gz}, \cite{fu}, \cite{fl} to write out the Hankel formula (LP--endpoints at $x=0$ and $x=\infty$ when $\nu\ge 1$). I.S. Kac in \cite{kac2} also  gave a general necessary and sufficient condition for a second order differential operator (including Sturm-Liouville operators, generalized second order derivatives of Krein-Feller type, and operators generated by canonical systems of phase dimension two) to have spectral multiplicity one, see also the more recent papers by D.J.Gilbert \cite{g} and B.Simon \cite{s}. 
\end{remark}

\begin{remark}
We also call attention to a recent paper \cite{vgt} of Voronov, Gitman and Tyutin on the Dirac equation with Coulomb potential. They also made use of M. Krein's method of directing functionals from \cite{k1} to obtain an eigenfunction expansion using a single solution of the Dirac equation.

\end{remark}

We finally mention that M. Krein was led to the method of directing functionals by some ideas in the dissertation of Moshe Liv\v sic from 1944, see \cite{liv}. In this dissertation also the {\it Stieltjes-Liv\v sic inversion formula}, see e.g. \cite{kk} and \cite[Lemma 3.2]{kst}, appeared for the first time.

\end{document}